\documentclass{article}
\pdfoutput=1
\usepackage{amsfonts}
\usepackage{geometry}
\usepackage{makeidx}
\usepackage{color}
\usepackage{amsthm,amsmath,amssymb}
\usepackage{mathrsfs}
\begin{document}
\title {MZV evaluation of hypergeometric series}
\author{Ming Hao Zhao}
\date{}

\maketitle
\begin{abstract}
We evaluate several classes of high weight hypergeometric series, with most of their parameters denominator-1/2/4, via Multiple Zeta Values of level 1/2/4.
\end{abstract}

\bibliographystyle{plain}

\tableofcontents

\section{Introduction}

We reduce certain hypergeometric series to 1/2/4-admissible polylog integrals (\textcolor{blue}{1APL/2APL/4APL}) or quadratic nested sums (\textcolor{blue}{QNS}), then evaluate them via colored Multiple Zeta Values (CMZVs, or simply MZVs) of level 1/2/4. Theoretical basis of this procedure, basic definitions of CMZV and admissibility, as well as irreducible constants' representation, is given by \cite{Au} and its associated algorithm. Here by QNS we mean all series of form $S=\sum_{R} \frac{(\pm1)^{n_1}\cdots (\pm1)^{n_k}}{f_1(n_1)^{s_1}\cdots f_k(n_k)^{s_k}}$, where $R$ is a restriction on indexes e.g. $n_1>\cdots>n_k>0$ (CMZVs)/$n_1=\max_j n_j$ (quadratic Euler sums)/etc, all $f_j(n) \in \{n+\mathbb{Z}, 2n+\mathbb{Z}\}$ and the weight $W=\sum_j s_j$. Evidently QNSs is reducible to level 4 CMZVs. For irreducible constants, we denote level 2/4 irreducible CMZVs as \textcolor{blue}{MZ, QMZ} respectively, say:\\
$$\text{MZ}(\{5,1\},\{-1,1\})=\sum_{n>m>0}\frac{(-1)^n}{n^5 m},\ \text{QMZ}(4,\{3,1,1\},\{0,0,1\})=\sum_{n>m>l>0}\frac{i^l}{n^3 m l}$$\\

\clearpage
\section{Argument 1}
\subsection{Denominator 2}
\noindent\textbf {Lemma 1.} Let $a,b \in \mathbb R$, $\ _{r+1}F_r\left(\{a\}_r, b; \{a+1\}_r;1\right)$ is expressible via polygamma functions. If $2a, 2b\in\mathbb Z$ (resp. $4a, 4b\in\mathbb Z$), then the series reduces to $\pi, \zeta(2k+1), \log(2)$ (resp. $\pi, \zeta(2k+1), \log(2), \beta(2j)$).\\

\noindent Proof. Evaluate $\int_0^1 t^p (1-t)^q \log ^n(t) \, dt$ in 2 ways: Power expansion on $(1-t)^b$/Beta derivatives. $\hfill\square$\\

\noindent\textbf {Lemma 2.} $\, _{r+1}F_r\left(\{1\}_{r+1};\frac{3}{2},\{2\}_{r-1};1\right)$ is expressible via level 2 MZVs.\\

\noindent Proof. By generalizing  \cite{ZMH} one have

$$\small \ _{r+3}F_{r+2}\left(\{1\}_{r+3};\frac32, \{2\}_{r+1};1\right)=\frac{2^{r+1} (-1)^r }{r!}\sum _{j=0}^r \binom{r}{j} (-\log (2))^{r-j}\int_0^{\frac{\pi }{2}} x \log ^j(2 \sin (x)) \, dx$$
Using contour integration, logsine integrals on RHS are reduced to 2APLs, then MZVs. \\

\noindent Proof 2. Expand $\text{Li}_n$ in the 2APL $\int_0^1 \frac{\text{Li}_n\left(\frac{4 x}{(x+1)^2}\right)}{2 x} \, dx$, integrate term by term.\\

\noindent Proof 3. Evaluate $\int_0^1 \frac{\text{Li}_n(t)}{t \sqrt{1-t}} \, dt$ in 2 ways: Power expansion on $\text{Li}_n$/Substitution $t\to 1-u^2$, 2APL. $\hfill\square$\\

\noindent\textbf {Lemma 3.} $ \, _{r+1}F_r\left(\{ \frac12\}_{r-1},\{1\}_2; \{ \frac32\}_r;1\right)$ is expressible via level 4 MZVs.\\

\noindent Proof 1. Evaluate $\int_0^1 \frac{t^a \sin ^{-1}(t) \log ^n(t)}{\sqrt{1-t^2}} \, dt$ in 2 ways: Power expansion on $\frac{\sin ^{-1}(t)}{\sqrt{1-t^2}}$/Substitution $t\to \frac{2u}{1+u^2}$, 4APL.\\

\noindent Proof 2. Apply double integration on $\int_0^1 \frac{\text{Li}_n(t)}{\sqrt{1-t^2}} \, dt$, (generalizing \cite{ZMH}). \\

\noindent Proof 3. Evaluate $\int_0^1 \frac{\cos ^{-1}(t) \text{Li}_n(t)}{t} \, dt$ in 2 ways: Power expansion on $\text{Li}_n$ (recall moment integral $\int_0^1 x^{m-1} \cos ^{-1}(x) \, dx=\frac{\sqrt{\pi } \Gamma \left(\frac{m+1}{2}\right)}{m^2 \Gamma \left(\frac{m}{2}\right)}$)/Substitution  $t\to \frac{2u}{1+u^2}$, 4APL.   $\hfill\square$\\

\noindent\textbf {Lemma 4.} $\pi  \, _{r+1}F_r\left(\{1\}_{r-1},\{ \frac32\}_2;\{2\}_r;1\right), \pi  \, _{r+1}F_r\left(\{ \frac12\}_{r+1};1,\{ \frac32\}_{r-1};1\right)$ is expressible via level 4 MZVs.\\

\noindent Proof:  The first class: By Power expansion on $K$ one have:\\
$$\ _{r+3}F_{r+2}\left( \{1\}_{r+1}, \frac32, \frac32; \{2\}_{r+2};1\right)=\frac{8(-1)^{r-1}}{\pi (r-1)!} \int_0^1 \frac{1}{t}\left(K(t)-\frac\pi 2\right) \log^ {r-1}(t) \, dt$$\\
And by FL theory one have:\\
$$\frac{K(x)-\frac{\pi }{2}}{x}=\sum _{n=0}^{\infty } (-1)^{n} (2 n+1) \left(4 \sum _{m=1}^n \frac{1}{m}\left(\frac{\pi }{4}-\sum _{k=0}^{m-1} \frac{(-1)^k}{2 k+1}\right)-4 C+2 \pi  \log (2)\right)P_n(2 x-1) $$\\
Notice its general coefficients are nested harmonic numbers. Since Theorem 3 of \cite{Jac} ensures that all FL coefficients of $\frac{\log^{r-1}(x)}x$ (hence $\log^r(x)$ by integration) are of the same structure, one may apply FL Parseval on $\log^r(x), \frac{K(x)-\frac{\pi }{2}}{x}$ to reduce the series to QNS, then MZVs.\\

The second class: FL Parseval on $K(x), \frac{\log^{r-1}(x)}{\sqrt{x}}$ (for FL of latter see \cite{Jac} again), QNS. $\hfill\square$\\

\noindent\textbf {Lemma 5.}  $_{r+1}F_r\left(\{1\}_{r+1};\{ \frac32\}_2,\{2\}_{r-2};1\right), \ \ _{r+1}F_r\left(\{\frac12\}_{r-2}, \{1\}_3;\{ \frac32\}_r;1\right)$ is expressible via level 4 MZVs.\\

\noindent Proof. The first class: Evaluate $\int_0^1 \frac{K(1-x) \text{Li}_n(x)}{x} \, dx$ in 2 ways: FL Parseval on $\frac{\text{Li}_n(x)}{x}, K(1-x)$, QNS/Power expansion on $\text{Li}_n(x)$ (recall moment integral $\int_0^1 K(1-x) x^{m-1} \, dx=\frac{\pi  \Gamma (m)^2}{2 \Gamma \left(m+\frac{1}{2}\right)^2}$).\\

The second class:  One may prove by induction. Case $k=2$ is trivial. Now assume it holds for $k=2,\cdots,r-1$. Using Beta integral:\\
$$ _{r+1}F_r\left(\{\frac12\}_{r-2}, \{1\}_3;\{ \frac32\}_r;1\right)\to \int_{0<x<1} \frac{dx}{\sqrt{x(1-x)}}\  \sqrt{x} \ _{r}F_{r-1}\left(\{\frac12\}_{r-2}, \{1\}_2;\{ \frac32\}_{r-1};x\right)$$\\
Here and below, by $A\to B$ we mean $A,B$ are equivalent modulo constants and trivial integrals/series. IBP:\\
$$\to \int_{0<x_1<1} \left(\int_{0<x<x_1}\frac{dx}{x_1\sqrt{x(1-x)}}\right)\ \sqrt{x_1}\  _{r-1}F_{r-2}\left(\{\frac12\}_{r-3}, \{1\}_2;\{ \frac32\}_{r-2};x_1\right)dx_1 $$\\
IBP induction:\\
$$\to \int_{0<x_{r-2}<1} \left(\int_{0<x<\cdots<x_{r-2}}\frac{dx_{r-3}\cdots dx_1 dx}{x_{r-2}\cdots x_1\sqrt{x(1-x)}}\right)\ \sqrt{x_{r-2}}\  _2F_1(1,1;\frac32;x_{r-2}) dx_{r-2}$$\\
Note that while calculating boundary values, closed-forms of case $k=2,\cdots,r-1$ are make of use. We also need that of $\int_{0<x<\cdots<x_{k}<1}\frac{dx_{k}\cdots dx_1 dx}{x_{k}\cdots x_1\sqrt{x(1-x)}}$ for $k=0,\cdots, r-3$, which is amenable using the same substitution below. Simplify $_2F_1$ and the integral:\\
$$\to \int_{0<x<\cdots<x_{r-2}<1}\frac{dx}{\sqrt{x(1-x)}}\left(\prod_{k-1}^{r-3} \frac{dx_k}{x_k}\right)\frac{\sin^{-1}(\sqrt{x_{r-2}})dx_{r-2}}{x_{r-2}\sqrt{1-x_{r-2}}}$$\\
Let $x\to \sin^2(t), x_k\to \sin^2(t_k)$:\\
$$\to \int_{0<t<\cdots<t_{r-2}<\frac{\pi}2}dt\left(\prod_{k-1}^{r-3} \cot(t_k)dt_k\right)\frac{t_{r-2}dt_{r-2}}{\sin(t_{r-2})}$$\\
Let $t\to 2\tan^{-1}(u), t_k\to 2\tan^{-1}(u_k)$:\\
$$\to \int_{0<u<\cdots<u_{r-2}<1}\frac{du}{1+u^2}\left(\prod_{k-1}^{r-3} \frac{(1-u_k^2)du_k}{u_k(1+u_k^2)}\right)\frac{\tan^{-1}(u_{r-2})du_{r-2}}{u_{r-2}}$$\\
The last integral is easily seen reducible to level 4 MZVs.  $\hfill\square$\\

\noindent \textbf{Theorem 1.}  Let $A$ (resp. $M, N, B$) be a vector with all components in $\mathbb Z/2$ (resp. $\mathbb N, \mathbb N, \mathbb C$), $A, M$ and $B, N$ are of same length, $S, T$ vectors that met one of five following conditions ($k,m,n,i,j\in\mathbb Z$):\\
\begin{itemize}
\item {$S=\{k\},\ T=\emptyset$}
\item {$ S=\{k+1/2\},\ T=\emptyset$}
\item {$S=\{k,m\},\ T=\{n+1/2\}$}
\item {$S=\{k+1/2, m+1/2\},\ T=\{n\}$}
\item {$S=\{k,m,n\},\ T=\{i+1/2,j+1/2\}$}\\
\end{itemize}

 \noindent Then $\, _{q+1}F_q(S,A,B;T,A+M,B-N;1)$, whenever convergent and non-terminating, is expressible via level 4 MZVs.\\

\noindent \textbf{Theorem 1'.} For all rational $R(n)$ satisfying all its poles belong to $\{-\frac12,-1,-\frac32,-2,\cdots \}\cup\{\frac12, \frac32,\frac52,\cdots\}$, all 5 classes of  binomial sums\\
$$\sum _{n=0}^{\infty } R(n) \left(\frac{\binom{2 n}{n}}{4^n}\right)^k\ (k=0, \pm1, -2),\ \pi \sum _{n=0}^{\infty } R(n) \left(\frac{\binom{2 n}{n}}{4^n}\right)^2$$\\
whevever convergent, are expressible via level 4 MZVs. \\

\noindent Proof. Using $a_n=\frac{\binom{2 n}{n}}{4^n}=\frac{\left(\frac{1}{2}\right)_n}{(1)_n}$ and partial fractions, equivalence between Thm 1 and 1' is justified. We only deal with 1' below.
Case $k=0$ corresponds to trivial rational series. Now suppose $k=1$, using $a_{n+1}=\frac{2n+1}{2n+2} a_n$, $a_{n-1}=\frac{2n}{2n-1}a_n$ one have\\
$$(1)\ \sum_n \frac{a_n}{2n+j}\to \sum_n \frac{2n a_n}{(2n-1)(2n+j-2)}\to \sum_n \frac{a_n}{2n-1}, \sum_n \frac{a_n}{2n+j-2}$$
$$(2)\ \sum_n \frac{a_n}{2n+j}\to \sum_n \frac{(2n+1) a_n}{(2n+2)(2n+j+2)}\to \sum_n  \frac{a_n}{2n+2}, \sum_n  \frac{a_n}{2n+j+2}$$\\
Here by $X\to Y$ we mean sums in $X,Y$ are equivalent modulo reindexing by $\pm1$, finite initial terms, partial fractions and trivial series. By induction using (1)(2), all $\sum_n \frac{a_n}{2n+j}$ are reduced to case $j=1,2$. In general, all $\sum_n \frac{a_n}{(2n+j)^m}$ are reduced to $\sum_n \frac{a_n}{(2n+1)^s}, \sum_n \frac{a_n}{(2n+2)^s},\ s=1,\cdots, m$, all of which are evaluable via Lemma 1. Finally by partial fractions one may evaluate all $\sum _{n=0}^{\infty } R(n) a_n$.\\

Now let $k=-1$. Since $a_n^{-1}=O(\sqrt{n})$, $\sum_n \frac{a_n^{-1}}{2n+j}$ is not convergent. We consider $\sum_n \frac{a_n^{-1}}{(2n+j)(2n+l)}$ instead. Using\\
$$\sum_n \frac{a_n^{-1}}{(2n+j)(2n+l)}\to \sum_n \frac{a_n^{-1}}{2n(2n+j-2)},  \sum_n \frac{a_n^{-1}}{2n(2n+l-2)}$$\\
The problem reduces to evaluation of $\sum_n \frac{a_n^{-1}}{2n(2n+j)}$. Reindex and do partial fractions appropriately:\\
$$(1')\ \sum_n \frac{a_n^{-1}}{2n(2n+j)}\to \sum_n \frac{a_n^{-1}}{2n(2n-2)}, \sum_n \frac{a_n^{-1}}{2n(2n+j-2)}$$
$$(2')\ \sum_n \frac{a_n^{-1}}{2n(2n+j)}\to \sum_n \frac{a_n^{-1}}{(2n+1)(2n+j+2)}$$ $$\to \sum_n \frac{a_n^{-1}}{(2n+1)(2n+j+4)}, \sum_n \frac{a_n^{-1}}{(2n+3)(2n+j+4)}\to \sum_n \frac{a_n^{-1}}{2n(2n+1)}, \sum_n \frac{a_n^{-1}}{2n(2n+j+2)}$$\\
Using (1')(2'), all $\sum_n \frac{a_n^{-1}}{2n(2n+j)}$ are reduced to case $j=1,2$. Similarly all $\sum_n \frac{a_n^{-1}}{(2n+j)^m}$ are reduced to $\sum_n \frac{a_n^{-1}}{(2n+1)^s}$, $\sum_n \frac{a_n^{-1}}{(2n+2)^s}, \ s=2,\cdots, m$.  Moreover, because of\\
$$\sum_n \frac{a_n^{-1}}{n^m}\to \sum_n \frac{a_n^{-1}}{(2n+2)^{m-1} (2n+1)}\to \sum_n \frac{a_n^{-1}}{(2n+1)(2n+2)}, \sum_n \frac{a_n^{-1}}{(2n+2)^{s}}\ (s=2,\cdots, m-1)$$\\ 
The following two classes $\sum_{n=1} \frac{a_n^{-1}}{n^m},\ \sum_{n=0} \frac{a_n^{-1}}{(2n+2)^m}$ can represent each other, thus it all boils down to $\sum_{n=0} \frac{a_n^{-1}}{(2n+1)^m},\ \sum_{n=1} \frac{a_n^{-1}}{n^m}$, which are solved in Lemma 2, 3. The general case follows by partial fractions again (one may combine several divergent $\sum_n \frac{a_n^{-1}}{2n+j}$ to convergent $\sum_n \frac{a_n^{-1}}{(2n+j)(2n+l)}$, which is possible assuming convergence).\\

Proof of case $k=2$ is similar to that of $k=1$, where in the last step one may use Lemma 4. For case $k=-2$, imitating $k=-1$, the basic units should be $\sum_n \frac{a_n^{-2}}{(2n+i)(2n+j)(2n+l)}$, $\sum_n \frac{a_n^{-2}}{(2n+j)^2(2n+l)}$ (due to convergence issues) and $\sum_n \frac{a_n^{-2}}{(2n+l)^m}\ (m\ge3)$, ultimately reducible to $\sum_{n=0} \frac{a_n^{-2}}{(2n+1)^m}$ and $\sum_{n=1} \frac{a_n^{-2}}{n^m}$ by repeated reindexing and partial fractions, evaluable using Lemma 5. $\hfill \square$\\

\subsection{Denominator 4}

\noindent \textbf{Theorem 2.} Let $a,b\in \mathbb Z/4, a\ge-1, b\ge0, r\in \mathbb N$ that met one of three following conditions:\\
\begin{itemize}
\item {$a\in \mathbb Z\ \text{or}\ b\in \mathbb {Z}$}
\item {$a+b\in \mathbb {Z}$}
\item {$a=b=-\frac14\ \text{or}\ -\frac34,\ r\le6$}\\
\end{itemize}

\noindent Then $_{r+1}F_r\left(\{1\}_{r},2+a;\{2\}_{r-1},3+a+b;1\right)$, whenever non-terminating and convergent, is expressible via level 4 MZVs.\\

\noindent Proof. We only need to evaluate $\int_0^1 x^a (1-x)^b \text{Li}_{r-1}(x) \, dx$. For $a\in \mathbb Z$ let $x\to1-t^4$, IBP, 4APL. For $b\in \mathbb Z$ let $x\to t^4$, IBP, 4APL. For $a+b\in \mathbb {Z}$, use partial fractions and Lemma 1. For the last class, apply FL method again; WLOG, assume $a=b=-\frac14, r=6$.\\

\noindent Step 1: FL expansions. Let $f(x)=\sum_{n=0}^\infty c_n P_n(2x-1)$, by solving difference equation of Legendre recurrence one have:\\
$$\int_0^x \frac{f(t)}t dt=\int_0^1 \frac{f(t)}t (1-t)\, dt+ \sum_{n=1}^\infty\left( (-1)^n \left(\frac1n+\frac1{n+1}\right) \sum_{k=n}^\infty(-1)^k c_k- \frac{c_n}{n+1}\right) P_n(2x-1)$$\\
Therefore, by using FL expansion of $\text{Li}_3(x)$ \cite{Jac}, one may derive that of $\text{Li}_n(x)$ by reindexing. We only need the following:\\

{\footnotesize$$\text{Li}_4(x)=\sum _{n=1}^{\infty } a_n P_n(2 x-1)-\zeta (3)+\frac{\pi ^4}{90}+\frac{\pi ^2}{6}-1$$\\

$$a_n=-2 (-1)^n \left(\frac{1}{n}+\frac{1}{n+1}\right) \sum _{k=n}^{\infty } \frac{(-1)^k}{k^3}-2 (-1)^n \left(\frac{1}{n^2}+\frac{1}{(n+1)^2}+\frac{2}{n}-\frac{2}{n+1}\right) \sum _{k=n}^{\infty } \frac{(-1)^k}{k^2}$$$$+4 (-1)^n \left(\frac{1}{n}+\frac{1}{n+1}\right) \sum _{k=n}^{\infty } \frac{1}{k}\sum _{j=k}^{\infty } \frac{(-1)^j}{j^2}+\frac{1}{n^4}+\frac{2}{n^3}-\frac{2}{n}+\frac{2}{n+1}+\frac{2}{(n+1)^2}-\frac{1}{(n+1)^4}$$}\\

{\footnotesize $$ \text{Li}_5(x)=\sum _{n=1}^{\infty } b_n P_n(2 x-1)+\zeta (3)+\zeta (5)-\frac{\pi ^2}{6}-\frac{\pi ^4}{90}+1$$\\

$$b_n=(-1)^n \left(-\frac{4}{n+1}-\frac{4}{n}\right) \sum _{k=n}^{\infty } \frac{1}{k}\sum _{j=k}^{\infty } \frac{(-1)^j}{j^3}+(-1)^n \left(-\frac{4}{n+1}-\frac{4}{n}\right) \sum _{k=n}^{\infty } \frac{1}{k^2}\sum _{j=k}^{\infty } \frac{(-1)^j}{j^2}$$$$+(-1)^n \left(-\frac{4}{n^2}+\frac{8}{n+1}-\frac{4}{(n+1)^2}-\frac{8}{n}\right) \sum _{k=n}^{\infty } \frac{1}{k}\sum _{j=k}^{\infty } \frac{(-1)^j}{j^2}+(-1)^n \left(\frac{8}{n+1}+\frac{8}{n}\right) \sum _{k=n}^{\infty } \frac{1}{k}\sum _{j=k}^{\infty } \frac{1}{j}\sum _{l=j}^{\infty } \frac{(-1)^l}{l^2}$$$$+(-1)^n \left(\frac{2}{n+1}+\frac{2}{n}\right) \sum _{k=n}^{\infty } \frac{(-1)^k}{k^4}+(-1)^n \left(\frac{2}{n^2}-\frac{4}{n+1}+\frac{2}{(n+1)^2}+\frac{4}{n}\right) \sum _{k=n}^{\infty } \frac{(-1)^k}{k^3}$$$$+(-1)^n \left(\frac{2}{n^3}+\frac{4}{n^2}-\frac{4}{(n+1)^2}+\frac{2}{(n+1)^3}\right) \sum _{k=n}^{\infty } \frac{(-1)^k}{k^2}-\frac{1}{n^5}-\frac{2}{n^4}-\frac{2}{n+1}-\frac{2}{(n+1)^2}-\frac{2}{(n+1)^3}+\frac{1}{(n+1)^5}+\frac{2}{n}$$}\\

\noindent Step 2: Integral representations. By repeated IBP one have:\\

{\footnotesize $$(-1)^n\sum _{k=n}^{\infty } \frac{(-1)^k}{k^2}=-\int_0^1 \frac{x^n \log (x)}{x (x+1)} \, dx,\ (-1)^n\sum _{k=n}^{\infty } \frac{(-1)^k}{k^3}=\frac{1}{2} \int_0^1 \frac{x^n \log ^2(x)}{x (x+1)} \, dx,\ (-1)^n\sum _{k=n}^{\infty } \frac{(-1)^k}{k^4}=-\frac{1}{6} \int_0^1 \frac{x^n \log ^3(x)}{x (x+1)} \, dx$$

$$(-1)^n\sum _{k=n}^{\infty } \frac{1}{k}\sum _{j=k}^{\infty } \frac{(-1)^j}{j^2}=\int_0^1 \frac{x^n }{x (x+1)}\left(-\text{Li}_2(-x)+\frac{\log ^2(x)}{2}-\log (x+1) \log (x)-\frac{\pi ^2}{12}\right) \, dx$$

$$(-1)^n\sum _{k=n}^{\infty } \frac{1}{k}\sum _{j=k}^{\infty } \frac{(-1)^j}{j^3}=\frac{1}{2} \int_0^1 \frac{x^n }{x (x+1)}\left(-2 \text{Li}_3(-x)+2 \text{Li}_2(-x) \log (x)-\frac{1}{3} \log ^3(x)+\log (x+1) \log ^2(x)-\frac{3 \zeta (3)}{2}\right) \, dx$$

$$(-1)^n\sum _{k=n}^{\infty } \frac{1}{k^2}\sum _{j=k}^{\infty } \frac{(-1)^j}{j^2}=-\int_0^1 \frac{x^n }{x (x+1)}\left(-2 \text{Li}_3(-x)+\text{Li}_2(-x) \log (x)+\frac{\log ^3(x)}{6}-\frac{1}{12} \pi ^2 \log (x)-\frac{3 \zeta (3)}{2}\right) \, dx$$

$$(-1)^n\sum _{k=n}^{\infty } \frac{1}{k}\sum _{j=k}^{\infty } \frac{1}{j}\sum _{l=j}^{\infty } \frac{(-1)^l}{l^2}=-\int_0^1 \frac{g(x) x^n}{x (x+1)} \, dx$$
$$g(x)=-\text{Li}_3(-x)-\text{Li}_3(x+1)+\frac{\log ^3(x)}{6}-\frac{1}{2} \log (x+1) \log ^2(x)-\frac{1}{2} i \pi  \log ^2(x+1)-\frac{1}{12} \pi ^2 \log (x)+\frac{1}{4} \pi ^2 \log (x+1)+\frac{\zeta (3)}{8}$$}\\

\noindent Step 3: Generating functions. By repeated integration one have:\\

{\footnotesize$$\sum _{n=1}^{\infty } \frac{\binom{2 n}{n} x^{2 n}}{4^n (2 n)}=\log \left(\frac{2}{\sqrt{1-x^2}+1}\right),\ \sum _{n=1}^{\infty } \frac{\binom{2 n}{n} x^{2 n}}{4^n (2 n+1)}=\frac{\sin ^{-1}(x)}{x}-1$$\\

$$\sum _{n=1}^{\infty } \frac{\binom{2 n}{n} x^{2 n}}{4^n (2 n+1)^2}=\frac{\log (2 x) \sin ^{-1}(x)}{x}+\frac{1}{4 i x}\left(\text{Li}_2\left(e^{2 i \sin ^{-1}(x)}\right)-\text{Li}_2\left(e^{-2 i \sin ^{-1}(x)}\right)\right)-1$$\\

$$\sum _{n=1}^{\infty } \frac{\binom{2 n}{n} x^{2 n}}{4^n (2 n)^2}=-\frac{1}{2} \text{Li}_2\left(\frac{1}{2} \left(\sqrt{1-x^2}+1\right)\right)-\frac{1}{2} \log \left(1-\sqrt{1-x^2}\right) \log \left(\frac{1}{2} \left(\sqrt{1-x^2}+1\right)\right)$$$$-\frac{1}{4} \log \left(\frac{1}{8} \left(\sqrt{1-x^2}+1\right)\right) \log \left(\frac{1}{2} \left(\sqrt{1-x^2}+1\right)\right)+\frac{\pi ^2}{12}$$}

{\scriptsize $$\sum _{n=1}^{\infty } \frac{\binom{2 n}{n} x^{2 n}}{4^n (2 n)^3}=\frac{1}{2} i \pi  \log (2) \log (x)+\frac{1}{12} \pi ^2 \log (x)+\frac{\zeta (3)}{4}+\frac{i \pi ^3}{24}+\frac{\log ^3(2)}{6}-\frac{1}{4} i \pi  \log ^2(2)-\frac{1}{24} \pi ^2 \log (2)$$$$-\frac{1}{4} i \pi  \text{Li}_2\left(\frac{1}{2} \left(1-\sqrt{1-x^2}\right)\right)-\frac{1}{4} i \pi  \text{Li}_2\left(\frac{1}{2} \left(\sqrt{1-x^2}+1\right)\right)+\frac{1}{4} \text{Li}_3\left(\frac{1}{2} \left(1-\sqrt{1-x^2}\right)\right)-\frac{1}{4} \text{Li}_3\left(\frac{1}{2} \left(\sqrt{1-x^2}+1\right)\right)$$$$-\frac{1}{2} \text{Li}_2\left(\frac{1}{2} \left(1-\sqrt{1-x^2}\right)\right) \log (x)+\frac{1}{4} \text{Li}_2\left(\frac{1}{2} \left(1-\sqrt{1-x^2}\right)\right) \log \left(\sqrt{1-x^2}+1\right)-\frac{1}{2} \text{Li}_2\left(\frac{1}{2} \left(\sqrt{1-x^2}+1\right)\right) \log (x)$$$$+\frac{1}{4} \text{Li}_2\left(\frac{1}{2} \left(\sqrt{1-x^2}+1\right)\right) \log \left(\sqrt{1-x^2}+1\right)-\frac{1}{24} \log ^3\left(\sqrt{1-x^2}+1\right)-\frac{3}{8} \log ^2(2) \log \left(1-\sqrt{1-x^2}\right)$$$$-\frac{3}{8} \log ^2(2) \log \left(\sqrt{1-x^2}+1\right)+\frac{1}{8} \log \left(1-\sqrt{1-x^2}\right) \log ^2\left(\sqrt{1-x^2}+1\right)+\frac{1}{2} \log (2) \log (x) \log \left(1-\sqrt{1-x^2}\right)$$$$+\frac{1}{2} \log (2) \log (x) \log \left(\sqrt{1-x^2}+1\right)-\frac{1}{2} \log (x) \log \left(1-\sqrt{1-x^2}\right) \log \left(\sqrt{1-x^2}+1\right)-\frac{1}{4} i \pi  \log \left(1-\sqrt{1-x^2}\right) \log \left(\sqrt{1-x^2}+1\right)$$}\\

\noindent Step 4. FL Parseval. By Dixon $_3F_2$ one have:\\
$$(x(1-x))^{s-1}=B(s,s)\sum_{n=0}^\infty  \frac{(5/4)_n(1-s)_n(1/2)_n}{(1/4)_n(1/2+s)_n(1)_n} P_{2n}(2x-1)$$\\
Therefore FL Parseval yields\\
$$\frac{1}{B\left(\frac{3}{4},\frac{3}{4}\right)}\int_0^1 \frac{\text{Li}_5(x)}{\sqrt[4]{x (1-x)}} \, dx=\sum _{n=1}^{\infty } \frac{\binom{2 n}{n} b_{2n}}{4^n (4 n+1)}-\zeta (2)+\zeta (3)-\zeta (4)+\zeta (5)+1$$\\

\noindent Step 5. Series-integral transformation. Observe that\\

{\scriptsize $$\frac{1}{4 n+1}\underset{n\to 2 n}{\text{lim}}\left(\frac{1}{n}+\frac{1}{n+1}\right)=\frac{1}{2 n}-\frac{1}{2 n+1},\ \frac{1}{4 n+1}\underset{n\to 2 n}{\text{lim}}\left(\frac{1}{(n+1)^2}+\frac{2}{n}+\frac{1}{n^2}-\frac{2}{n+1}\right)=\frac{1}{4 n^2}-\frac{1}{(2 n+1)^2}$$

$$\frac{1}{4 n+1}\underset{n\to 2 n}{\text{lim}}\left(\frac{1}{(n+1)^3}+\frac{2}{n^2}+\frac{1}{n^3}-\frac{2}{(n+1)^2}\right)=\frac{1}{8 n^3}-\frac{1}{(2 n+1)^3}$$}

\noindent And that all coefficients of nested sums in expression of $b_n$ are constant times of those in brackets, $\sum _{n=1}^{\infty } \frac{\binom{2 n}{n} b_{2n}}{4^n (4 n+1)}$ is reducible to a combination of following modulo trivial rational series:\\
$$\sum _{n=1}^{\infty } \frac{\binom{2 n}{n} }{4^n (2 n)^k}\int_0^1 \frac{f(x) x^{2n}}{x (x+1)} \, dx,\ \sum _{n=1}^{\infty } \frac{\binom{2 n}{n} }{4^n (2 n+1)^k}\int_0^1 \frac{f(x) x^{2n}}{x (x+1)} \, dx$$\\
Where $k=1,2,3$ and $f(x)$ one of the 7 functions in Step 2 above. Apply Fubini based on 5 generating functions in Step 3, all components except one are transformed into integrals with polylog-arcsine integrands. The only exception $k=1, f(x)=-\log(x)$, which is not directly transformable since no polylog closed-form is known for $\sum _{n=1}^{\infty } \frac{\binom{2 n}{n} x^{2 n}}{4^n (2 n+1)^3}$, also reduces to the previous case after IBP:\\

{\footnotesize $$\frac{-1}{(2n+1)^3}\int_0^1 \frac{x^{2 n+1} \log (x)}{(x+1) x^2} \, dx=\frac{1}{(2 n+1)^2}\int_0^1  x^{2 n} \left(\text{Li}_2(-x)-\frac{1}{x}-\frac{1}{2} \log ^2(x)+\log (x+1) \log (x)-\frac{\log (x)}{x}+\frac{\pi ^2}{12}+1\right)\, dx$$}\\

\noindent Step 6. 4-admissible transform. Substitution $x\to \sin(t), t\to 2\tan^{-1}(u)$ transforms the resulting integral in last step to a nonhomogeneous 4APL, thus 4APLs after IBP, finally MZVs. Case $r\le5$ and $a=b=-\frac34, r=6$ are easier. This completes the proof. $\hfill\square$\\

\section{Other results}

\noindent \textbf{Proposition 1.}  Let $A$ (resp. $A_1, A_2, M, N, B$) be a vector with all components in $\mathbb Z/2$ (resp. $\mathbb Z/4, \mathbb R, \mathbb N, \mathbb N, \mathbb C$), $A, A_1,A_2,  M$ and $B, N$ are of same length, $k\in\mathbb N$, $t\in \mathbb R$.  Then the hypergeometric series $\, _{q+1}F_q\left(k,A,B;A+M,B-N;-1\right)$ and $\, _{q+1}F_q\left(k,A_1,B;A_1+M,B-N;1\right)$ (resp. generalized $\, _{q+1}F_q\left(t,A_2,B;A_2+M,B-N;1\right)$), whenever convergent and non-terminating, is expressible via level 4 MZVs (resp. polygamma functions).\\

\noindent Proof. First 2 cases correspond to trivial rational series $\sum_{n=0}^{\infty} \frac{(-1)^n}{\prod_{i=1}^N (2n+a_i)^{m_i}}$ and $\sum_{n=0}^{\infty} \frac{1}{\prod_{i=1}^N (4n+a_i)^{m_i}}$. The last case is also direct by combining result of Lemma 1 and argument of partial fractions given in Thm 1. Note that in the last step i.e. evaluating sth like $_3F_2(a,b,c;b+1,c+1;1)$ ($1<a<2$), one may replace the argument 1 with $z$, do partial fractions, let $z\to 1^-$, use Abel theorem and asymptotics of $ _2F_1$ to circumvent the divergence issue. $\hfill\square$\\

\noindent \textbf{Proposition 2.}  Let $A$ (resp. $M, N, B$) be a vector with all components in $X_A$ (resp. $\mathbb N, \mathbb N, \mathbb C$), $A, M$ and $B, N$ are of same length, $S, T$ vectors that met one of three following conditions ($k,m,n\in\mathbb Z$):\\
\begin{itemize}
\item {$S=\{k\},\ T=\emptyset,\ X_A=\mathbb Z$}
\item {$ S=\{k+1/2\},\ T=\emptyset,\ X_A=\mathbb Z+\frac12$}
\item {$S=\{k,m\},\ T=\{n+1/2\},\ X_A=\mathbb Z$}\\
\end{itemize}

 \noindent Then $\, _{q+1}F_q\left(S,A,B;T,A+M,B-N;\frac12\right)$, whenever convergent and non-terminating, is expressible via level 4 MZVs (in fact level 2 for the first case).\\

\noindent Proof. The first case corresponds to trivial polylog series. For last 2 cases, start by evaluating $\int_0^1 t^p \sin^{-1}(\sqrt{t/2})^k (1-\frac{t}2)^{-\frac12} \log ^n(t) \, dt$ in 2 ways for $k=0,1$ respectively: Power expansion on $\sin^{-1}(\sqrt{t/2})^k (1-\frac{t}2)^{-\frac12}$/Substution $t\to \frac{2u^2}{1+u^2}$. The rest are similar to Thm. 1. $\hfill\square$\\

\noindent \textbf{Proposition 3.} $\, _{r+1}F_r\left(\{1\}_{r+1};\frac{3}{2},\{2\}_{r-1};-\frac18\right)$ is expressible via level 2 MZVs.\\

\noindent Proof. By generalizing \cite{ZMH} one have
$$\small \ _{r+1}F_r\left(\{1\}_{r+1};\frac32, \{2\}_{r-1};-\frac18\right)=\frac{2 (-1)^{r-1} }{(r-3)!}\int_0^1 \frac{(y+2) \log ^2(y+1) \log ^{r-3}\left(\frac{2 y^2}{y+1}\right)}{y (y+1)} \, dy$$\\
Note that 2APLs on RHS are reducible to  MZVs.$\hfill\square$\\

\noindent \textbf{Proposition 4.} $\, _{r+1}F_r\left(\{1\}_r,\frac{3}{2};\frac{4}{3},\frac{5}{3},\{2\}_{r-2};\frac{2}{27}\right)$ is expressible via level 4 MZVs.\\

\noindent Proof. Since our original solution is ad hoc and limited, we refer to \cite{Au} for a systematic proof. $\hfill \square$\\

We conclude by 2 conjectures:\\

\noindent \textbf{Conjecture 1.}  Let $A, B$ be vectors with all components in $\mathbb N$ and of same length, $k,m,n\in\mathbb Z$, then $_{r+1}F_r\left(k,A,\frac{m}4;B,\frac{n}4;1\right)$ is expressible via level 4 MZVs.\\

\noindent \textbf{Conjecture 2.}  Let $S$ be a $_{q+1}F_q$ with all its parameters in $\mathbb Z/2$ and its closed-form expressible via rational combination of level 4 MZVs, then $S$ belong to one of the five classes in Thm. 1.\\\\

\section{Examples}

Lemma 1:

{\footnotesize $$\, _{10}F_9\left(\{1\}_9,\frac{3}{2};\{2\}_9;1\right)=\frac{2 \pi ^2 \zeta (3)^2}{3}-24 \zeta (3) \zeta (5)-\frac{16}{15} \zeta (3) \log ^5(2)+\frac{8}{9} \pi ^2 \zeta (3) \log ^3(2)-16 \zeta (5) \log ^3(2)-8 \zeta (3)^2 \log ^2(2)$$ $$+\frac{1}{5} \pi ^4 \zeta (3) \log (2)+4 \pi ^2 \zeta (5) \log (2)-72 \zeta (7) \log (2)+\frac{2339 \pi ^8}{907200}-\frac{4}{315}  \log ^8(2)+\frac{4}{135} \pi ^2 \log ^6(2)+\frac{1}{30} \pi ^4 \log ^4(2)+\frac{79 \pi ^6 \log ^2(2)}{3780}$$\\

$$\, _{10}F_9\left(\{ \frac12\}_{10}; \{ \frac32\}_9;1\right)=\frac{\pi ^3 \zeta (3)^2}{1536}+\frac{3 \pi  \zeta (3) \zeta (5)}{128}+\frac{1}{960} \pi  \zeta (3) \log ^5(2)+\frac{\pi ^3 \zeta (3) \log ^3(2)}{1152}+\frac{1}{64} \pi  \zeta (5) \log ^3(2)+\frac{1}{128} \pi  \zeta (3)^2 \log ^2(2)$$ $$+\frac{19 \pi ^5 \zeta (3) \log (2)}{46080}+\frac{1}{256} \pi ^3 \zeta (5) \log (2)+\frac{9}{128} \pi  \zeta (7) \log (2)+\frac{11813 \pi ^9}{928972800}+\frac{\pi  \log ^8(2)}{80640}+\frac{\pi ^3 \log ^6(2)}{34560}+\frac{19 \pi ^5 \log ^4(2)}{276480}+\frac{55 \pi ^7 \log ^2(2)}{774144}$$\\

$$\frac{1024 \sqrt{\pi } }{\Gamma \left(\frac{1}{4}\right)^2}\, _7F_6\left(\{\frac14\}_6,\frac{3}{4};\{\frac54\}_6;1\right)=7 C \zeta (3)+\frac{23 \pi ^3 C}{48}+\pi  C^2+\frac{1}{6} C \log ^3(2)+\frac{1}{4} \pi  C \log ^2(2)+2 C^2 \log (2)+\frac{3}{8} \pi ^2 C \log (2)$$ $$+\frac{21 \pi ^2 \zeta (3)}{32}+\frac{93 \zeta (5)}{4}+\frac{7}{8} \zeta (3) \log ^2(2)+\frac{7}{8} \pi  \zeta (3) \log (2)+\frac{1}{32} \pi  \zeta \left(4,\frac{1}{4}\right)+\frac{1}{16} \zeta \left(4,\frac{1}{4}\right) \log (2)+\frac{587 \pi ^5}{5120}+\frac{\log ^5(2)}{480}+\frac{1}{192} \pi  \log ^4(2)$$ $$+\frac{1}{64} \pi ^2 \log ^3(2)+\frac{23}{384} \pi ^3 \log ^2(2)+\frac{73 \pi ^4 \log (2)}{1536}$$\\

$$\sqrt{\frac{2}{\pi ^3}} \Gamma \left(\frac{1}{4}\right)^2 \, _7F_6\left(\frac{1}{4},\{\frac12\}_6;\{\frac32\}_6;1\right)=\frac{7 C \zeta (3)}{2}+\frac{11 \pi ^3 C}{96}+4 C^2+\frac{\pi ^2 C}{4}-\frac{\pi  C^2}{2}+4 \pi  C-32 C+\frac{1}{12} C \log ^3(2)+\frac{1}{8} \pi  C \log ^2(2)$$ $$-C \log ^2(2)-C^2 \log (2)-\frac{1}{16} \pi ^2 C \log (2)-\pi  C \log (2)+8 C \log (2)+\frac{7 \pi ^2 \zeta (3)}{64}+\frac{7 \pi  \zeta (3)}{4}-14 \zeta (3)-\frac{93 \zeta (5)}{8}-\frac{7}{16} \zeta (3) \log ^2(2)$$ $$+\frac{7}{2} \zeta (3) \log (2)-\frac{7}{16} \pi  \zeta (3) \log (2)-\frac{\zeta \left(4,\frac{1}{4}\right)}{8}+\frac{1}{64} \pi  \zeta \left(4,\frac{1}{4}\right)+\frac{1}{32} \zeta \left(4,\frac{1}{4}\right) \log (2)+\frac{27 \pi ^4}{256}-\pi ^2-\frac{11 \pi ^3}{24}-\frac{1241 \pi ^5}{30720}-16 \pi +128$$ $$-\frac{1}{960} \log ^5(2)+\frac{\log ^4(2)}{48}-\frac{1}{384} \pi  \log ^4(2)+\frac{1}{384} \pi ^2 \log ^3(2)-\frac{\log ^3(2)}{3}+\frac{1}{24} \pi  \log ^3(2)-\frac{11}{768} \pi ^3 \log ^2(2)-\frac{1}{32} \pi ^2 \log ^2(2)-\frac{1}{2} \pi  \log ^2(2)$$ $$+4 \log ^2(2)+\frac{11}{96} \pi ^3 \log (2)+\frac{1}{4} \pi ^2 \log (2)-\frac{27 \pi ^4 \log (2)}{1024}+4 \pi  \log (2)-32 \log (2)$$}\\

\clearpage

\noindent Lemma 2:

{\footnotesize $$\, _9F_8\left(\{1\}_9;\frac{3}{2},\{2\}_7;1\right)=-\frac{4}{9} \pi ^2 \text{MZ}(\{5,1\},\{-1,1\})-\frac{26}{3} \text{MZ}(\{7,1\},\{-1,1\})-\frac{8}{3} \text{MZ}(\{5,1,1,1\},\{-1,1,-1,1\})$$ $$-\frac{32}{3} \text{Li}_5\left(\frac{1}{2}\right) \zeta (3)-\frac{2}{27} \pi ^4 \text{Li}_4\left(\frac{1}{2}\right)+64 \text{Li}_8\left(\frac{1}{2}\right)+\frac{\pi ^2 \zeta (3)^2}{3}+\frac{251 \zeta (3) \zeta (5)}{16}+\frac{4}{45} \zeta (3) \log ^5(2)+\frac{14}{27} \pi ^2 \zeta (3) \log ^3(2)+\frac{31}{36} \zeta (5) \log ^3(2)$$ $$-\frac{53}{540} \pi ^4 \zeta (3) \log (2)+\frac{247}{72} \pi ^2 \zeta (5) \log (2)+\frac{1651}{96} \zeta (7) \log (2)-\frac{76357 \pi ^8}{10886400}+\frac{\log ^8(2)}{630}+\frac{2}{135} \pi ^2 \log ^6(2)-\frac{67 \pi ^4 \log ^4(2)}{3240}-\frac{853 \pi ^6 \log ^2(2)}{45360}$$}\\

\noindent Lemma 3:

{\footnotesize $$ \, _6F_5\left(\{ \frac12\}_4,\{1\}_2; \{ \frac32\}_5;1\right)=2 \Im(\text{QMZ}(4,\{4,1\},\{1,0\}))-2 \Im(\text{QMZ}(4,\{4,1\},\{1,2\}))$$ $$+16 \Im\left(\text{Li}_5\left(\frac{1}{2}+\frac{i}{2}\right)\right)-\frac{35 \pi ^5}{1536}-\frac{1}{96} \pi  \log ^4(2)-\frac{1}{64} \pi ^3 \log ^2(2)$$\\}

\noindent Lemma 4:

{\footnotesize $$\pi  \, _7F_6\left(\{1\}_5,\{ \frac32\}_2;\{2\}_6;1\right)=-2560 \Im(\text{QMZ}(4,\{4,1\},\{1,0\}))+\frac{9728}{3} \Im(\text{QMZ}(4,\{4,1\},\{1,2\}))-16384 \Im\left(\text{Li}_5\left(\frac{1}{2}+\frac{i}{2}\right)\right)$$ $$-64 \pi  \zeta (3) \log (2)+\frac{4}{3} \zeta \left(4,\frac{1}{4}\right) \log (2)-\frac{4}{3} \zeta \left(4,\frac{3}{4}\right) \log (2)+25 \pi ^5-32 \pi  \log ^4(2)+48 \pi ^3 \log ^2(2)$$\\

$$\pi  \, _6F_5\left(\{ \frac12\}_6;1,\{ \frac32\}_4;1\right)=-40 \Im(\text{QMZ}(4,\{4,1\},\{1,0\}))+\frac{152}{3} \Im(\text{QMZ}(4,\{4,1\},\{1,2\}))-256 \Im\left(\text{Li}_5\left(\frac{1}{2}+\frac{i}{2}\right)\right)$$ $$+\frac{1}{48} \zeta \left(4,\frac{1}{4}\right) \log (2)-\frac{1}{48} \zeta \left(4,\frac{3}{4}\right) \log (2)+\frac{25 \pi ^5}{64}+\frac{1}{6} \pi  \log ^4(2)+\frac{3}{4} \pi ^3 \log ^2(2)$$}\\

\noindent Lemma 5:

{\footnotesize $$\, _6F_5\left(\{1\}_6;\{ \frac32\}_2,\{2\}_3;1\right)=128 \pi  \Im\left(\text{Li}_4\left(\frac{1}{2}+\frac{i}{2}\right)\right)-64 \text{Li}_5\left(\frac{1}{2}\right)+\frac{217 \zeta (5)}{4}$$ $$-\frac{3}{8} \pi  \zeta \left(4,\frac{1}{4}\right)+\frac{3}{8} \pi  \zeta \left(4,\frac{3}{4}\right)+\frac{8 \log ^5(2)}{15}-\frac{2}{9} \pi ^2 \log ^3(2)+\frac{41}{45} \pi ^4 \log (2)$$\\

$$\, _6F_5\left(\{\frac12\}_3, \{1\}_3;\{ \frac32\}_5;1\right)=-16 \pi  \Im\left(\text{Li}_4\left(\frac{1}{2}+\frac{i}{2}\right)\right)+16 \text{Li}_5\left(\frac{1}{2}\right)-\frac{341 \zeta (5)}{32}$$ $$+\frac{3}{64} \pi  \zeta \left(4,\frac{1}{4}\right)-\frac{3}{64} \pi  \zeta \left(4,\frac{3}{4}\right)-\frac{2}{15}  \log ^5(2)+\frac{5}{36} \pi ^2 \log ^3(2)-\frac{37}{360} \pi ^4 \log (2)$$}\\

\clearpage

\noindent Theorem 1 (one for each case):\\

{\footnotesize $$\, _7F_6\left(\frac{1}{2},1,\{\frac54\}_5;\frac{3}{2},\{\frac94\}_5;1\right)=-\frac{3125 C}{81}-\frac{96875 \zeta (5)}{96}-\frac{21875 \zeta (3)}{216}-\frac{3125 \zeta \left(4,\frac{1}{4}\right)}{1152}$$ $$+\frac{756250}{243}-\frac{3125 \pi ^2}{648}-\frac{3125 \pi ^3}{864}-\frac{3125 \pi }{972}-\frac{15625 \pi ^5}{4608}-\frac{1}{486} 3125 \log (2)$$ \\

$$\, _8F_7\left(\{\frac12\}_4,\frac{7}{6},\frac{5}{4},\frac{4}{3},\frac{3}{2};\frac{1}{6},\frac{1}{4},\frac{1}{3},\{\frac52\}_4;1\right)=\frac{2835 \pi  \zeta (3)}{32}-\frac{17739 \pi }{128}-\frac{1593 \pi ^3}{512}$$ $$+\frac{945}{16} \pi  \log ^3(2)-\frac{4779}{128} \pi  \log ^2(2)+\frac{945}{64} \pi ^3 \log (2)-\frac{3645}{64} \pi  \log (2)$$\\

$$\, _8F_7\left(\{\frac12\}_4,1,1,\frac{4}{3},\frac{5}{3};\frac{1}{3},\frac{2}{3},\{\frac32\}_4,\frac{5}{2};1\right)=-\frac{3}{8} \Im(\text{QMZ}(4,\{4,1\},\{1,0\}))+\frac{3}{8} \Im(\text{QMZ}(4,\{4,1\},\{1,2\}))$$ $$-\frac{105 C}{64}+\frac{105}{16} \Im\left(\text{Li}_3\left(\frac{1}{2}+\frac{i}{2}\right)\right)+\frac{3}{4} \Im\left(\text{Li}_4\left(\frac{1}{2}+\frac{i}{2}\right)\right)-3 \Im\left(\text{Li}_5\left(\frac{1}{2}+\frac{i}{2}\right)\right)+\frac{3 \zeta \left(4,\frac{3}{4}\right)}{2048}-\frac{3 \zeta \left(4,\frac{1}{4}\right)}{2048}$$ $$+\frac{35 \pi ^5}{8192}+\frac{105}{128}-\frac{105 \pi ^3}{2048}+\frac{1}{512} \pi  \log ^4(2)+\frac{1}{256} \pi  \log ^3(2)+\frac{3 \pi ^3 \log ^2(2)}{1024}-\frac{105}{512} \pi  \log ^2(2)+\frac{3 \pi ^3 \log (2)}{1024}$$\\

$$\pi  \, _7F_6\left(\{-\frac12\}_2,\{1\}_5;\{2\}_6;1\right)=-\frac{2560}{9} \Im(\text{QMZ}(4,\{4,1\},\{1,0\}))+\frac{9728}{27} \Im(\text{QMZ}(4,\{4,1\},\{1,2\}))$$ $$-\frac{47104 C}{243}-\frac{14336}{27} \Im\left(\text{Li}_3\left(\frac{1}{2}+\frac{i}{2}\right)\right)-\frac{32768}{27} \Im\left(\text{Li}_4\left(\frac{1}{2}+\frac{i}{2}\right)\right)-\frac{16384}{9} \Im\left(\text{Li}_5\left(\frac{1}{2}+\frac{i}{2}\right)\right)+\frac{256 \pi  \zeta (3)}{27}$$ $$-\frac{64}{9} \pi  \zeta (3) \log (2)+\frac{32 \zeta \left(4,\frac{1}{4}\right)}{9}-\frac{32 \zeta \left(4,\frac{3}{4}\right)}{9}+\frac{4}{27} \zeta \left(4,\frac{1}{4}\right) \log (2)-\frac{4}{27} \zeta \left(4,\frac{3}{4}\right) \log (2)+\frac{25 \pi ^5}{9}+\frac{112 \pi ^3}{9}$$ $$-\frac{46784 \pi }{729}+\frac{117248}{729}-\frac{32}{9}  \pi  \log ^4(2)+\frac{512}{27} \pi  \log ^3(2)+\frac{16}{3} \pi ^3 \log ^2(2)-\frac{448}{9} \pi  \log ^2(2)-\frac{128}{9} \pi ^3 \log (2)+\frac{23552}{243} \pi  \log (2)$$\\

$$\, _7F_6\left(\{1\}_6,\frac{3}{2};\{2\}_3,\{\frac52\}_3;1\right)=1512 \pi  C+2592 \pi  \Im\left(\text{Li}_3\left(\frac{1}{2}+\frac{i}{2}\right)\right)+3456 \pi  \Im\left(\text{Li}_4\left(\frac{1}{2}+\frac{i}{2}\right)\right)$$ $$-2592 \text{Li}_4\left(\frac{1}{2}\right)-1728 \text{Li}_5\left(\frac{1}{2}\right)-3024 \zeta (3)+\frac{5859 \zeta (5)}{4}-\frac{81}{8} \pi  \zeta \left(4,\frac{1}{4}\right)+\frac{81}{8} \pi  \zeta \left(4,\frac{3}{4}\right)-\frac{369 \pi ^4}{10}$$ $$-1620 \pi +4536+\frac{72 \log ^5(2)}{5}-108 \log ^4(2)-6 \pi ^2 \log ^3(2)+27 \pi ^2 \log ^2(2)+\frac{123}{5} \pi ^4 \log (2)$$}

\clearpage

\noindent Theorem 2 (one for each case):\\

{\footnotesize $$\, _6F_5\left(\{1\}_6;\{2\}_4,\frac{11}{4};1\right)=392 \Im(\text{QMZ}(4,\{4,1\},\{1,0\}))-\frac{1456}{3} \Im(\text{QMZ}(4,\{4,1\},\{1,2\}))-336 \Re(\text{QMZ}(4,\{3,1,1\},\{0,0,1\}))$$ $$+336 C \Im\left(\text{Li}_3\left(\frac{1}{2}+\frac{i}{2}\right)\right)-98 C \zeta (3)+\frac{112 C^2}{3}-\frac{28 \pi ^2 C}{3}-\frac{49 \pi ^3 C}{8}+\frac{896 C}{27}-\frac{21}{2} \pi  C \log ^2(2)+21 \pi ^2 C \log (2)+\frac{448}{3} \Im\left(\text{Li}_3\left(\frac{1}{2}+\frac{i}{2}\right)\right)$$ $$+672 \Im\left(\text{Li}_4\left(\frac{1}{2}+\frac{i}{2}\right)\right)+3024 \Im\left(\text{Li}_5\left(\frac{1}{2}+\frac{i}{2}\right)\right)-112 \text{Li}_4\left(\frac{1}{2}\right)-147 \text{Li}_5\left(\frac{1}{2}\right)+105 \text{Li}_4\left(\frac{1}{2}\right) \log (2)+\frac{1211 \pi ^2 \zeta (3)}{64}+\frac{11151 \zeta (5)}{128}$$ $$-\frac{392 \zeta (3)}{9}-\frac{77 \zeta \left(4,\frac{1}{4}\right)}{48}-\frac{21}{64} \pi  \zeta \left(4,\frac{1}{4}\right)+\frac{77 \zeta \left(4,\frac{3}{4}\right)}{48}+\frac{21}{64} \pi  \zeta \left(4,\frac{3}{4}\right)-\frac{35}{192} \zeta \left(4,\frac{1}{4}\right) \log (2)+\frac{35}{192} \zeta \left(4,\frac{3}{4}\right) \log (2)+\frac{119 \pi ^4}{240}-\frac{49 \pi ^3}{18}$$ $$-\frac{112 \pi ^2}{27}+\frac{1792}{243}-\frac{2415 \pi ^5}{512}+\frac{28 \log ^5(2)}{5}-\frac{14 \log ^4(2)}{3}-\frac{63}{32} \pi  \log ^4(2)+\frac{35}{24} \pi ^2 \log ^3(2)+\frac{7}{2} \pi  \log ^3(2)-\frac{441}{64} \pi ^3 \log ^2(2)-\frac{35}{6} \pi ^2 \log ^2(2)$$ $$-\frac{14}{3} \pi  \log ^2(2)-\frac{21}{160} \pi ^4 \log (2)+\frac{49}{8} \pi ^3 \log (2)+\frac{28}{3} \pi ^2 \log (2)$$\\

$$\, _7F_6\left(\{1\}_6,\frac{9}{4};\{2\}_5,3;1\right)=\frac{288 C \zeta (3)}{5}+\frac{28 \pi ^2 C}{5}-\frac{192 C^2}{5}-\frac{32 \pi  C^2}{5}-\frac{96 C}{5}-\frac{14 \pi ^3 C}{15}+\frac{48 \pi  C}{5}+\frac{144}{5} C \log ^3(2)$$ $$-\frac{72}{5} \pi  C \log ^2(2)-\frac{432}{5} C \log ^2(2)+\frac{192}{5} C^2 \log (2)-\frac{28}{5} \pi ^2 C \log (2)+\frac{144}{5} \pi  C \log (2)-\frac{288}{5} C \log (2)-\frac{21 \pi ^2 \zeta (3)}{5}+\frac{108 \pi  \zeta (3)}{5}$$ $$+\frac{792 \zeta (5)}{5}-\frac{216 \zeta (3)}{5}+\frac{324}{5} \zeta (3) \log ^2(2)-\frac{648}{5} \zeta (3) \log (2)-\frac{108}{5} \pi  \zeta (3) \log (2)-\frac{3 \zeta \left(4,\frac{1}{4}\right)}{5}-\frac{1}{10} \pi  \zeta \left(4,\frac{1}{4}\right)+\frac{3 \zeta \left(4,\frac{3}{4}\right)}{5}$$ $$+\frac{1}{10} \pi  \zeta \left(4,\frac{3}{4}\right)+\frac{3}{5} \zeta \left(4,\frac{1}{4}\right) \log (2)-\frac{3}{5} \zeta \left(4,\frac{3}{4}\right) \log (2)+\frac{243 \pi ^4}{400}+\frac{7 \pi ^3}{10}+\frac{7 \pi ^2}{5}+\frac{8}{5}-\frac{1219 \pi ^5}{7200}+\frac{12 \pi }{5}+\frac{81 \log ^5(2)}{25}-\frac{81 \log ^4(2)}{5}$$ $$-\frac{27}{10} \pi  \log ^4(2)-\frac{21}{10} \pi ^2 \log ^3(2)-\frac{108 \log ^3(2)}{5}+\frac{54}{5} \pi  \log ^3(2)-\frac{21}{20} \pi ^3 \log ^2(2)+\frac{63}{10} \pi ^2 \log ^2(2)-\frac{108 \log ^2(2)}{5}+\frac{54}{5} \pi  \log ^2(2)$$ $$-\frac{243}{400} \pi ^4 \log (2)+\frac{21}{10} \pi ^3 \log (2)+\frac{21}{5} \pi ^2 \log (2)-\frac{72 \log (2)}{5}+\frac{36}{5} \pi  \log (2)$$\\

$$\, _7F_6\left(\{1\}_6,\frac{7}{4};\{2\}_5,\frac52;1\right)=36 \Im(\text{QMZ}(4,\{4,1\},\{1,0\}))-20 \Im(\text{QMZ}(4,\{4,1\},\{1,2\}))-64 \Re(\text{QMZ}(4,\{3,1,1\},\{0,0,1\}))$$ $$+64 C \Im\left(\text{Li}_3\left(\frac{1}{2}+\frac{i}{2}\right)\right)+21 C \zeta (3)-2 \pi ^3 C+\frac{10 \pi ^2 C}{3}-8 \pi  C^2+16 \pi  C-\frac{4}{3} C \log ^3(2)-8 C \log ^2(2)+16 C^2 \log (2)$$ $$+\frac{5}{3} \pi ^2 C \log (2)+8 \pi  C \log (2)-32 C \log (2)-64 \Im\left(\text{Li}_3\left(\frac{1}{2}+\frac{i}{2}\right)\right)+64 \Im\left(\text{Li}_4\left(\frac{1}{2}+\frac{i}{2}\right)\right)-64 \Im\left(\text{Li}_5\left(\frac{1}{2}+\frac{i}{2}\right)\right)$$ $$+44 \text{Li}_4\left(\frac{1}{2}\right)-2 \text{Li}_5\left(\frac{1}{2}\right)+20 \text{Li}_4\left(\frac{1}{2}\right) \log (2)-\frac{37 \pi ^2 \zeta (3)}{16}-14 \pi  \zeta (3)-21 \zeta (3)-\frac{457 \zeta (5)}{64}+7 \zeta (3) \log ^2(2)-7 \pi  \zeta (3) \log (2)$$ $$+28 \zeta (3) \log (2)-\frac{\zeta \left(4,\frac{1}{4}\right)}{8}+\frac{1}{16} \pi  \zeta \left(4,\frac{1}{4}\right)+\frac{\zeta \left(4,\frac{3}{4}\right)}{8}-\frac{1}{16} \pi  \zeta \left(4,\frac{3}{4}\right)-\frac{7}{32} \zeta \left(4,\frac{1}{4}\right) \log (2)+\frac{7}{32} \zeta \left(4,\frac{3}{4}\right) \log (2)$$ $$+\frac{95 \pi ^5}{384}+2 \pi ^3-\frac{10 \pi ^2}{3}-\frac{277 \pi ^4}{480}-16 \pi +64+\frac{13 \log ^5(2)}{15}+2 \log ^4(2)-\frac{67}{72} \pi ^2 \log ^3(2)+\frac{4 \log ^3(2)}{3}+\frac{1}{4} \pi ^3 \log ^2(2)$$ $$-\frac{9}{4} \pi ^2 \log ^2(2)+8 \log ^2(2)-\frac{97}{960} \pi ^4 \log (2)+\pi ^3 \log (2)-\frac{5}{3} \pi ^2 \log (2)-8 \pi  \log (2)+32 \log (2)$$}

\clearpage

\noindent Proposition 1 (one for each case):\\

{\footnotesize $$\, _9F_8\left(2,\{\frac12\}_4,\frac{1}{3},\frac{2}{3},i,-i;\frac{3}{2},\frac{3}{2},\frac{5}{2},\frac{5}{2},-\frac{1}{3},-\frac{2}{3},i-1,-i-1;-1\right)$$ $$=\frac{11619 C}{512}+\frac{4095 \zeta \left(4,\frac{1}{4}\right)}{65536}-\frac{5373}{256}-\frac{3051 \pi ^3}{8192}-\frac{1365 \pi ^4}{16384}+\frac{6309 \pi }{4096}$$\\

$$\, _9F_8\left(1,-\frac{1}{2},-\frac{1}{2},-\frac{1}{4},-\frac{1}{4},\frac{1}{6},\frac{1}{6},i,i;\frac{3}{2},\frac{3}{2},\frac{3}{4},\frac{3}{4},-\frac{5}{6},-\frac{5}{6},i-1,i-1;1\right)$$ $$=-\left(\frac{49}{150}-\frac{343 i}{3600}\right) C+\left(\frac{19}{25}-\frac{83 i}{1350}\right)+\left(\frac{17}{300}-\frac{79 i}{9600}\right) \pi ^2-\left(\frac{7}{600}-\frac{343 i}{10800}\right) \pi +\left(\frac{7}{300}-\frac{343 i}{5400}\right) \log (2)$$\\

$$\, _8F_7\left(\frac{1}{2},\frac{1}{2},\frac{3}{4},\frac{3}{4},1,1,\frac{5}{4},\frac{5}{4};\frac{1}{4},\frac{3}{2},\frac{3}{2},\frac{7}{4},\frac{7}{4},2,3;1\right)$$ $$=\frac{156824}{75}+108 \pi -648 \log (2)+\frac{8 \sqrt{\pi ^3}}{\Gamma \left(\frac{1}{4}\right)^2} \left(-100 \sqrt{2}-\frac{288 \pi }{5}-\frac{8064}{25}-3 \sqrt{2} \pi -6 \sqrt{2} \log (2)+\frac{576 \log (2)}{5}\right)$$}

\noindent Proposition 2 (one for each case):\\

{\footnotesize $$\, _9F_8\left(\{1\}_9;\{2\}_4,\{3\}_4;\frac{1}{2}\right)=1184 \text{Li}_4\left(\frac{1}{2}\right)-640 \text{Li}_5\left(\frac{1}{2}\right)+320 \text{Li}_6\left(\frac{1}{2}\right)-128 \text{Li}_7\left(\frac{1}{2}\right)$$ $$+32 \text{Li}_8\left(\frac{1}{2}\right)-1120 \zeta (3)+416 \pi ^2-5280-\frac{1}{3} 640 \log ^3(2)-2496 \log ^2(2)+\frac{320}{3} \pi ^2 \log (2)+3840 \log (2)$$\\

$$\sqrt{2} \, _6F_5\left(-\frac{1}{2},\{\frac12\}_3,\frac{3}{2},\frac{3}{2};\{\frac52\}_5;\frac{1}{2}\right)=\frac{243}{128} \Im(\text{QMZ}(4,\{4,1\},\{1,0\}))-\frac{18711 C}{2048}+\frac{1215}{512} \Im\left(\text{Li}_3\left(\frac{1}{2}+\frac{i}{2}\right)\right)$$ $$-\frac{1215}{128} \Im\left(\text{Li}_4\left(\frac{1}{2}+\frac{i}{2}\right)\right)+\frac{243}{32} \Im\left(\text{Li}_5\left(\frac{1}{2}+\frac{i}{2}\right)\right)-\frac{1215 \pi  \zeta (3)}{1024}-\frac{243}{512} \pi  \zeta (3) \log (2)+\frac{1215 \zeta \left(4,\frac{1}{4}\right)}{131072}$$ $$-\frac{1215 \zeta \left(4,\frac{3}{4}\right)}{131072}+\frac{243 \zeta \left(4,\frac{1}{4}\right) \log (2)}{65536}-\frac{243 \zeta \left(4,\frac{3}{4}\right) \log (2)}{65536}+\frac{93555}{4096}-\frac{9315 \pi ^3}{65536}-\frac{56511 \pi ^5}{1310720}+\frac{101817 \pi }{16384}$$ $$-\frac{243 \pi  \log ^4(2)}{16384}-\frac{1215 \pi  \log ^3(2)}{8192}-\frac{3645 \pi  \log ^2(2)}{16384}-\frac{1863 \pi ^3 \log ^2(2)}{32768}-\frac{18711 \pi  \log (2)}{8192}-\frac{9315 \pi ^3 \log (2)}{32768}$$\\

$$\, _6F_5\left(\{1\}_6;-\frac{1}{2},2,2,3,3;\frac{1}{2}\right)=-176 \pi  C-112 \pi  \Im\left(\text{Li}_3\left(\frac{1}{2}+\frac{i}{2}\right)\right)-48 \pi  \Im\left(\text{Li}_4\left(\frac{1}{2}+\frac{i}{2}\right)\right)+140 \text{Li}_4\left(\frac{1}{2}\right)$$ $$+30 \text{Li}_5\left(\frac{1}{2}\right)-3 \pi ^2 \zeta (3)+385 \zeta (3)+\frac{1209 \zeta (5)}{16}+\frac{3}{64} \pi  \zeta \left(4,\frac{1}{4}\right)-\frac{3}{64} \pi  \zeta \left(4,\frac{3}{4}\right)+\frac{133 \pi ^4}{72}+\frac{15 \pi ^2}{2}+118 \pi $$ $$-296-\frac{1}{4} \log ^5(2)+\frac{35 \log ^4(2)}{6}-\frac{1}{12} \pi ^2 \log ^3(2)+\frac{7}{6} \pi ^2 \log ^2(2)-\frac{19}{48} \pi ^4 \log (2)-22 \pi ^2 \log (2)$$}\\

\clearpage

\noindent Proposition 3:\\

{\footnotesize $$\, _9F_8\left(\{1\}_9;\frac{3}{2},\{2\}_7;-\frac{1}{8}\right)=\frac{20}{9} \pi ^2 \text{MZ}(\{5,1\},\{-1,1\})-\frac{32}{3} \text{MZ}(\{7,1\},\{-1,1\})+\frac{40}{3} \text{MZ}(\{5,1,1,1\},\{-1,1,-1,1\})$$ $$-24 \log ^2(2) \text{MZ}(\{5,1\},\{-1,1\})-24 \log (2) \text{MZ}(\{5,1,1\},\{-1,1,1\})+\frac{160}{3} \text{Li}_5\left(\frac{1}{2}\right) \zeta (3)+\frac{10}{27} \pi ^4 \text{Li}_4\left(\frac{1}{2}\right)+112 \text{Li}_8\left(\frac{1}{2}\right)$$ $$+24 \text{Li}_7\left(\frac{1}{2}\right) \log (2)+\frac{5 \pi ^2 \zeta (3)^2}{6}-\frac{1351 \zeta (3) \zeta (5)}{16}-\frac{23}{45} \zeta (3) \log ^5(2)+\frac{20}{27} \pi ^2 \zeta (3) \log ^3(2)+\frac{269}{18} \zeta (5) \log ^3(2)-8 \zeta (3)^2 \log ^2(2)$$ $$+\frac{136}{135} \pi ^4 \zeta (3) \log (2)+\frac{133}{72} \pi ^2 \zeta (5) \log (2)+\frac{415}{6} \zeta (7) \log (2)-\frac{4499 \pi ^8}{340200}-\frac{19 \log ^8(2)}{10080}+\frac{1}{270} \pi ^2 \log ^6(2)+\frac{29 \pi ^4 \log ^4(2)}{3240}-\frac{103 \pi ^6 \log ^2(2)}{1134}$$}\\

\noindent Proposition 4:\\

{\footnotesize $$\, _7F_6\left(\{1\}_6,\frac{3}{2};\frac{4}{3},\frac{5}{3},\{2\}_4;\frac{2}{27}\right)=24 \pi  \Im\left(\text{Li}_4\left(\frac{1}{2}+\frac{i}{2}\right)\right)-153 \text{Li}_5\left(\frac{1}{2}\right)-90 \text{Li}_4\left(\frac{1}{2}\right) \log (2)+\frac{3 \pi ^2 \zeta (3)}{2}$$ $$+27 \zeta (5)-18 \zeta (3) \log ^2(2)+\frac{9}{128} \pi  \zeta \left(4,\frac{1}{4}\right)-\frac{9}{128} \pi  \zeta \left(4,\frac{3}{4}\right)-\frac{1}{40} 97 \log ^5(2)+\frac{41}{24} \pi ^2 \log ^3(2)-\frac{61}{160} \pi ^4 \log (2)$$}\\\\

\noindent \textbf {\large Acknowledgements.} Special thanks to K. C. Au, Chen and J. D'Aurizio for valuable suggestions, leading to improvements of this article.\\

\end{document}